\documentclass{comnet}

\usepackage{amssymb,amsfonts,amsmath,amsthm}
\usepackage{graphicx}
\usepackage{xcolor}
\usepackage[caption=false]{subfig}
\usepackage{url}
\usepackage{float}

\newcommand{\pd}[2]{\dfrac{\partial #1}{\partial #2}}
\newcommand{\ER}{Erd\H{o}s-R\'{e}nyi}
\newcommand{\BA}{Barab{\'a}si-Albert}
\newcommand{\x}{\mathbf{x}}
\renewcommand{\u}{\mathbf{u}}

\begin{document}

\title{The spectral degree exponent of a graph}

\shorttitle{Spectral degree exponent of a graph} 
\shortauthorlist{Achterberg and Van Mieghem} 

\author{
\name{Massimo A. Achterberg$^*$ and Piet Van Mieghem}
\address{Faculty of Electrical Engineering, Mathematics and Computer Science, Delft University of Technology, P.O. Box 5031, 2600 GA Delft, The Netherlands}
\email{$^*$Corresponding author: M.A.Achterberg@tudelft.nl}}

\maketitle

\begin{abstract}
{We propose the spectral degree exponent as a novel graph metric. Although Hofmeister \cite{HofmeisterThesis} has studied the same metric, we generalise Hofmeister's work to weighted graphs. We provide efficient iterative formulas and bounds for the spectral degree exponent and provide highly accurate asymptotic expansions for the spectral degree exponent for several families of graphs. Furthermore, we uncover a close relation between the spectral degree exponent and the well-known degree assortativity, by showing high correlations between the two metrics in all small graphs, several random graph models and many real-world graphs.}
{Graph metric, Assortativity, Similarity, Spectral radius, Degree sequence}
\end{abstract}

\section{Introduction}\label{sec_introduction}

We introduce the spectral degree exponent and show its close relation to the degree assortativity for undirected, but weighted and possibly disconnected networks. Similar to the degree assortativity, the spectral degree exponent also gives an indication whether nodes with similar degrees are connected to each other. Section~\ref{sec_introduce} formally introduces the spectral degree exponent (SDE) and discusses its historic background. Section~\ref{sec_extremal_values} investigates the extremal values of the SDE and its relation to the degree assortativity. Section~\ref{sec_numerics} discusses an iterative method to determine the SDE and provides bounds on the SDE. Section~\ref{sec_special_cases} computes the SDE explicitly for several graphs. Finally, Section~\ref{sec_simulations} demonstrates the similarity of the degree assortativity and spectral degree exponent using simulations on small graphs, \ER{}, \BA{} and many real-world communication networks. We conclude in Section~\ref{sec_conclusion}.

\section{The spectral degree exponent}\label{sec_introduce}
Let $G$ be a simple, undirected graph on $N$ nodes and with $L$ non-negatively weighted links. Let $d_{\max} = d_1 \geq d_2 \geq \ldots \geq d_N$ be the ordered weighted degree sequence of the graph $G$. The spectral radius~$\lambda_1$ is the largest eigenvalue of the weighted adjacency matrix $A$ of the graph $G$. We define the spectral degree exponent as

\begin{definition}\label{def_q}
	The \textbf{spectral degree exponent} (SDE) of a simple, undirected, non-regular (but possibly disconnected and weighted) graph $G$ is defined as the solution for $q$ of the equation
	\begin{equation}\label{eq_def}
		\lambda_1 = \left( \frac{1}{N} \sum_{i=1}^N d_i^{\, q} \right)^{1/q}.
	\end{equation}
\end{definition}
We say that the spectral degree exponent (SDE) is $q$ or has value $q$. Equation~\eqref{eq_def} can be rewritten as
\begin{equation}\label{eq_relation_dmax}
	\lambda_1 = \left( \frac{c}{N} d_{\max}^{\, q} + \frac{1}{N} \sum_{i=c+1}^N d_i^{\, q} \right)^{1/q}
\end{equation}
where $c$ is the multiplicity of the weighted degree $d_{\max}$. If $D$ denotes the random variable of the degree of node in the graph $G$, then a probabilistic rephrasing of definition (\ref{eq_def}) is
\begin{align}
	\lambda_{1}^{q}  &  =\sum_{k=d_{\min}}^{d_{\max}}\Pr\left[  D=k\right]
	k^{q} \label{def_SDE_probabilitistic}
\end{align}
The existence of $q$ in (\ref{def_q}) is established by the following, basic theorem:
\begin{theorem}[Existence and uniqueness \cite{HofmeisterThesis} {\cite[art.\ 74]{GS2}} ]\label{thm_q}
	The spectral degree exponent (SDE) exists, is unique and is bounded in $q \in [2,\infty]$.
\end{theorem}

Theorem~\ref{thm_q} has already been established for \emph{unweighted} graphs by Michael Hofmeister in his 1986 PhD thesis \cite{HofmeisterThesis}. Hofmeister originally named $q$ the ``Spektral Mittelwertcharakteristik'' or ``spectral mean characteristic'' due to its relation to the spectral radius and the q-norm of the degree sequence~\cite{hofmeister1988spectralradius}. Hofmeister \cite{HofmeisterThesis,hofmeister1988spectralradius} extended the notion of $q$ to graphs with multi-links and directed links. In his paper, we aim to provide an interpretation of the SDE $q$ in relation to the degree assortativity $\rho_D$, which will be motivated in Section~\ref{sec_extremal_values}. Up to our knowledge, the SDE $q$ only has a relevant interpretation for undirected, simple graphs. Hence, in the remainder of this work, we assume that the considered graphs are simple and undirected.

\section{Extremal values and the degree assortativity}\label{sec_extremal_values}
The SDE $q$ cannot be undefined for regular graphs. Indeed, for a regular graph $G$ with weighted degree~$r$, the definition~\eqref{eq_def} reduces to
\begin{align*}
	\lambda_1 &= \left( \frac{1}{N} \sum_{i=1}^N d_i^{\, q} \right)^{1/q} = \left( \frac{1}{N} \sum_{i=1}^N r^{\, q} \right)^{1/q} = \left( r^{\, q} \right)^{1/q} = r,
\end{align*}
which is independent of $q$ and, hence, the SDE is not defined. Therefore, regular graphs are excluded in Theorem~\ref{thm_q}.

Prior to continuing the analysis, we recall the well-known degree assortativity~$\rho_D$, which describes the Pearson correlation coefficient of nodes with different degrees. If all nodes have the same degree, the degree assortativity is undefined. In fact, as we will see later, the extremal values of the SDE and the degree assortativity are very similar.

The degree assortativity $\rho_D$ is maximally assortative (+1) when almost all nodes connect to nodes of the same degree and is maximally disassortative (-1) when all nodes connect to nodes of opposite degree. Prior to examining the values values $q=2$ and $q=\infty$ of the SDE, we introduce a certain graph family.

\begin{definition}
	A graph $G$ is \textbf{biregular} (also known as semiregular bipartite) if it is bipartite and the nodes in the same group have the same (weighted) degree.
\end{definition}
The complete bipartite graphs and the star graph are special cases of biregular graphs.

\begin{theorem}[Theorem 4.1 in \cite{hofmeister1988spectralradius}] Let $G$ be a simple, non-regular graph. \\
	\phantom{a}\qquad (a) If $G$ is biregular, then $q=2$. \\
	\phantom{a}\qquad (b) If $G$ is connected and $q=2$, then $G$ is biregular.
	\label{thm_q=2}
\end{theorem}
\textit{Proof.} See Appendix~\ref{app_proof_q2}.

\begin{theorem}[Satz II.9 in \cite{HofmeisterThesis}]
	Let $G$ be a simple, non-regular graph with maximum weighted degree $d_{\max}$. It holds that $q=\infty$ if and only if $G$ is disconnected and one component of $G$ is a clique of weighted size $d_{\max}+1$.
	\label{thm_q=infty}
\end{theorem}
\textit{Proof.} See Appendix~\ref{app_proof_qinfty}.

For the degree assortativity $\rho_D$, the most disassortative network is the complete bipartite graph $K_{m,n}$ where $m \neq n$. Theorem~\ref{thm_q=2} states that the minimal $q=2$ is attained in biregular graphs, of which the complete bipartite graphs are a special case. Similarly, the maximal $q=\infty$ for disconnected graphs is attained for a graph that has one regular component equal to the maximum weighted degree, whereas the most assortative graphs are regular graphs. Thus the extremal graphs for the degree assortativity and the SDE are very similar.

Finally, we discuss the concept of degree-preserving rewiring. In degree-preserving rewiring, two chosen links are rewired, such that the weighted degree of all nodes in the network remain the same. Van Mieghem \emph{et al}.\ \cite{vanmieghem2010degree} showed that applying degree-preserving rewiring to increase the degree assortativity~$\rho_D$ also increases the spectral radius~$\lambda_1$. Following the definition of the SDE $q$ in Eq.~\eqref{eq_def}, keeping the degrees fixed but increasing the spectral radius $\lambda_1$ directly leads to the increase of the SDE $q$. 

We conclude that the extremal cases as well as the important property of degree-preserving rewiring coincide well for the degree assortativity and the SDE. We further confirm the tight relationship between the degree assortativity and the SDE in Section~\ref{sec_simulations}.

\section{Characteristics of the spectral degree exponent}\label{sec_numerics}
Since the SDE $q$ in Equation \eqref{eq_def} is defined implicitly, the SDE cannot be computed explicitly for many graphs. Instead, this section focuses on numerical approaches, upper and lower bounds and a recursive formula.

\subsection{Root-finding}
Rewriting Definition \eqref{eq_def} leads to
\begin{equation}\label{eq_f}
	f(q) = \lambda_1 - \left( \frac{1}{N} \sum_{i=1}^N d_i^{\, q} \right)^{1/q}
\end{equation}
where $f(2) \geq 0$ and $f(\infty) \leq 0$ for any non-regular graph, we obtain the SDE $q$ as the root of $f(q)$. Unfortunately, finding the root of \eqref{eq_f} is numerically unstable, because $d_i^{\, q}$ can be extremely large. To guarantee stability, we take the logarithm of \eqref{eq_def} to find
\begin{equation*}
	q \log(\lambda_1) = -\log\left( N\right) + \log\left( \sum_{i=1}^N d_i^{\, q} \right)
\end{equation*}
and then define
\begin{equation}\label{eq_f2}
	f_1(q) = q \log(\lambda_1) + \log\left( N\right) - \log\left( \sum_{i=1}^N d_i^{\, q} \right)
\end{equation}
for which $f_1(2) \geq 0$ and $f_1(\infty) \leq 0$. The root of Eq.~\eqref{eq_f2} can be computed using the bisection method in a numerically stable way using \texttt{numpy.logaddexp}.

\subsection{Recursion}
We derive a recursive formula for the SDE~$q$ by rewriting Eq.~\eqref{eq_def} as
\[
N\lambda_{1}^{q}=\sum_{i=1}^{N}d_{i}^{\,q}=d_{\max}^{q}\left(  c+\sum
_{i=c+1}^{N}\left(  \frac{d_i}{d_{\max}}\right)
^{q}\right)
\]
such that
\begin{equation}\label{rewriting_def_SDE}
	\left( \frac{\lambda_{1}}{d_{\max}}\right)^{q} = \frac{1}{N}\left(c + \sum_{i=c+1}^{N}\left( \frac{d_i}{d_{\max}}\right)^{q}\right)
\end{equation}
We know that $\frac{\lambda_{1}}{d_{\max}}\leq1$ (from Gershgorin's
Theorem \cite[art.\ 245]{GS2}) and that $q\geq2$.
After taking the logarithm of both sides,%
\begin{equation}\label{q_log_transform}
	q = \frac{\log N-\log\left(  c + \sum_{i=c+1}^{N}\left(  \frac{d_i}{d_{\max}}\right)^{q}\right)  }{\log\left( \frac{d_{\max}}{\lambda_{1}}\right)}
\end{equation}
we can define the recursive formula
\begin{equation}\label{eq_q_recursion}
	q_k = \frac{\log N-\log\left(  c + \sum_{i=c+1}^{N}\left(  \frac
		{d_i}{d_{\max}}\right)^{q_{k-1}}\right)  }{\log\left(
		\frac{d_{\max}}{\lambda_{1}}\right)  }
\end{equation}
with initial starting value $q_{0}=\frac{\log \frac{N}{c}}{\log\left( \frac{d_{\max}}{\lambda_{1}}\right) }$, which is an upper bound for $q=\lim_{k\rightarrow\infty}q_{k}$. Since all $\frac{d_i}{d_{\max}}<1$ in the sum and $\left(  \frac{d_{\left(  i\right)  }}{d_{\max
}}\right)  ^{q_{k-1}}$ is a decreasing function in $q_{k-1}>0$, the recursion \eqref{eq_q_recursion} converges quickly in most cases. For example, for \ER{} graphs with $N=50$, $p=0.2$ and initial guess~$q_0$, at most 20 steps are required to obtain an accuracy of $10^{-6}$.

\subsection{Bounds}\label{sec_bounds}
Using the resursive formula \eqref{eq_q_recursion}, upper and lower bounds for the SDE can be obtained.
\begin{theorem}\label{Theorem_bounds}
	The spectral degree exponent (SDE) $q$ is bounded by
	\begin{equation}\label{eq_q_upper_bound2}
		\frac{\log N - \log\left( c + ( N-c) \left( \frac{d_{2}}{d_{\max}}\right)  ^{2} \right) }{\log\left( \frac{d_{\max}}{\lambda_{1}}\right)} \leq q \leq \frac{\log \frac{N}{c}}{\log\left( \frac{d_{\max}}{\lambda_{1}}\right) }
	\end{equation}
	where $c$ is the multiplicity of the maximum weighted degree $d_{\max}$ and $d_2$ is the second-largest weighted degree.
\end{theorem}
\textit{Proof.} The sum on the right-hand side of (\ref{rewriting_def_SDE}) is bounded as%
\[
\left(  N-c\right)  \left(  \frac{d_{\min}}{d_{\max}}\right)  ^{q}\leq
\sum_{i=c+1}^{N}\left(  \frac{d_{i}}{d_{\max}}\right)  ^{q}\leq\left(
N-c\right)  \left(  \frac{d_{2}}{d_{\max}}\right)  ^{q}%
\]
Since $\frac{d_{i}}{d_{\max}}<1$ for any $i>c$, the lower bound is further lower bounded for $q\rightarrow \infty$, while the upper bound is upper bounded for $q=2$, because the SDE has a value of at least $q=2$. Hence,
\[
0\leq\left(  N-c\right)  \left(  \frac{d_{\min}}{d_{\max}}\right)  ^{q}\leq
\sum_{i=c+1}^{N}\left(  \frac{d_{i}}{d_{\max}}\right)  ^{q}\leq \left(
N-c\right)  \left(  \frac{d_{2}}{d_{\max}}\right)  ^{q} \leq \left(
N-c\right)  \left(  \frac{d_{2}}{d_{\max}}\right)  ^{2}
\]
Substituting the most left and most right bound into (\ref{q_log_transform}) demonstrates Theorem \ref{Theorem_bounds}.  
\hfill $\square$

For graphs without isolated nodes, the proof can be sharpened. From the proof of Theorem \ref{Theorem_bounds} follows that for a graph with minimal weighted degree $d_{\min}>0$ the upper bound can be sharpened as
\begin{equation}\label{eq_q_upper_bound_better}
	q \leq \frac{\log N - \log\left( c + ( N-c) \left( \frac{d_{min}}{d_{\max}}\right)^{q_0} \right) }{\log\left( \frac{d_{\max}}{\lambda_{1}}\right)} < \frac{\log \frac{N}{c}}{\log\left( \frac{d_{\max}}{\lambda_{1}}\right) }
\end{equation}
since
\(
0<\left(  N-c\right)  \left(  \frac{d_{\min}}{d_{\max}}\right)  ^{q_0}\leq \left(  N-c\right)  \left(  \frac{d_{\min}}{d_{\max}}\right)  ^{q}
\).

\subsection{SDE is not monotonic}
Starting with a given network and adding a randomly chosen links usually increases the SDE $q$. The monotonicity of the SDE $q$ is, however, not guaranteed, as shown by the next example.

\begin{figure}[!ht]
	\centering
	\includegraphics[width=0.6\textwidth]{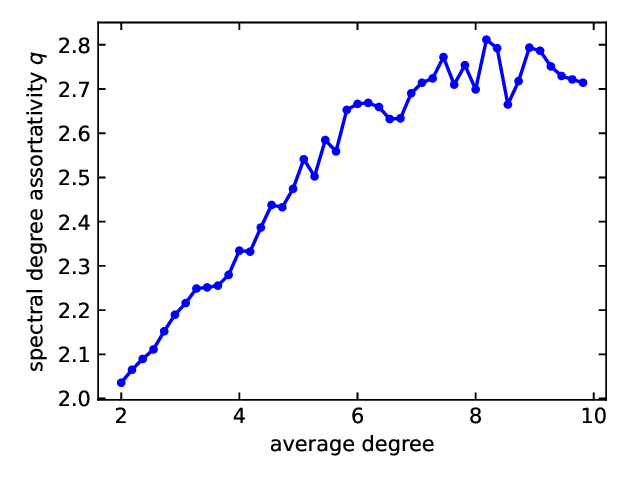}
	\caption{The SDE $q$ is not monotonic in adding links, as exemplified in a star graph on 11 nodes. Sequentially adding non-existing links at random increases $q$ for the majority of the time, but not in all cases.}
	\label{fig_q_non_monotonic}
\end{figure}

The star graph $G$ on $N$ nodes has the lowest possible SDE $q=2$ (see Theorem~\ref{thm_q=2}). Also, the maximum degree $d_{\max}=N-1$ in the star graph. We randomly add non-existing links to $G$ until $G$ equals the complete graph~$K_N$. Figure~\ref{fig_q_non_monotonic} shows that randomly adding links in a star graph with 11 nodes nearly always increases $q$, but not under all circumstances.

\section{The spectral degree exponent of some graphs}\label{sec_special_cases}
The spectral degree exponent~$q$ is computed (asymptotically) for several graph families.

\subsection{The path graph}
\begin{figure}[!ht]
	\centering
	\includegraphics{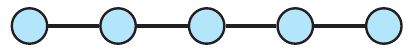}
	\caption{The path graph $P_N$. In this example, $N=5$.}
	\label{fig_path_graph}
\end{figure}
\begin{theorem}
	For the path graph $P_N$ on $N$ nodes, it holds that
	\begin{equation}\label{eq_q_path}
		q =\frac{4}{\pi^2} N + \frac{12}{\pi^2} + \left( \frac{1}{3} + \frac{52}{3\pi^{2}} \right) \frac{1}{N} + O\left(\frac{1}{N^2}\right), \qquad \textnormal{for } N \to\infty.
	\end{equation}
\end{theorem}
\textit{Proof.} For the path graph \cite{BrouwerHaemers}, the spectral radius $\lambda_1 = 2 \cos\left(\frac{\pi}{N+1}\right)$. Substituting $d_1 = \ldots = d_{N-2} = 2$ and $d_{N-1} = 1, d_N = 1$ into Eq.\ \eqref{eq_def}, we find
\begin{align}\label{eq_mfsyfosydfoys}
	\cos^q\left(\frac{\pi}{N+1}\right) &= 1 - \frac{2}{N} + \frac{2^{1-q}}{N}.
\end{align}
The Taylor expansion of the cosine for large $N$ is given by
\begin{equation*}
	\cos\left(\frac{\pi}{N+1}\right) = 1 - \frac{\pi^2}{2N^2} + \frac{\pi^2}{N^3} + \frac{\pi^2(\pi^2-36)}{24N^4} + O\left(\frac{1}{N^5}\right)
\end{equation*}
As an ansatz for large $N$, we use $q = a N + b + c/N$, where $a>0$ and $b,c \in \mathbb{R}$ are three parameters. Substituting back in \eqref{eq_mfsyfosydfoys} yields
\begin{equation*}
	\left( 1 - \frac{\pi^2}{2N^2} + \frac{\pi^2}{N^3} + \frac{\pi^2(\pi^2-36)}{24N^4} + O\left(\frac{1}{N^5}\right) \right)^{aN+b+c/N} = 1 - \frac{2}{N} + \frac{2^{1-aN+b+c/N}}{N}
\end{equation*}
For large $N$, the third term on the right-hand side is extremely small and can be neglected. We approximate the left-hand side using the binomial coefficient $(1+x)^N \approx 1 + xN + \frac{1}{2} x^2 N^2 + O(N^3)$, such that
\begin{align*}
	&1 + \left(aN+b+\frac{c}{N} \right) \left( \frac{-\pi^2}{2N^2} + \frac{\pi^2}{N^3} + \frac{\pi^2(\pi^2-36)}{24N^4} \right) + \frac{1}{2} \left(aN+b+\frac{c}{N} \right)^2 \left( \frac{-\pi^2}{2N^2} + \frac{\pi^2}{N^3} + \frac{\pi^2(\pi^2-36)}{24N^4} \right)^2 \\
	&\qquad + \frac{1}{6} \left(aN+b+\frac{c}{N} \right)^3 \left( \frac{-\pi^2}{2N^2} + \frac{\pi^2}{N^3} + \frac{\pi^2(\pi^2-36)}{24N^4} \right)^3 + O\left(\frac{1}{N^4}\right) = 1 - \frac{2}{N}
\end{align*}
Simplifying gives
\begin{align*}
	&1 - \frac{a \pi^2}{2} \frac{1}{N} + \left( -\frac{b \pi^2}{2} + a \pi^2 + \frac{a^2 \pi^4}{8} \right) \frac{1}{N^2} \\
	&\qquad + \left( -\frac{c \pi^2}{2} + b \pi^2 + \frac{a\pi^2(\pi^2-36)}{24} + \frac{ab\pi^4}{4} - \frac{a^2 \pi^4}{2} - \frac{a^3 \pi^6}{48} \right) \frac{1}{N^3} + O\left(\frac{1}{N^4}\right) = 1 - \frac{2}{N}
\end{align*}
Equating corresponding powers $k$ in $N^{-k}$ yields
\begin{align*}
	k=0: &\qquad 1 = 1 \\
	k=1: &\qquad -\frac{a \pi^2}{2} = -2 \\
	k=2: &\qquad -\frac{b \pi^2}{2} + a \pi^2 + \frac{a^2 \pi^4}{8} = 0 \\
	k=3: &\qquad -\frac{c \pi^2}{2} + b \pi^2 + \frac{a\pi^2(\pi^2-36)}{24} + \frac{ab\pi^4}{4} - \frac{a^2 \pi^4}{2} - \frac{a^3 \pi^6}{48} = 0
\end{align*}
which results in $a = \frac{4}{\pi^2}, b= \frac{12}{\pi^2}$ and $c=\frac{1}{3} + \frac{52}{3\pi^{2}}$.\hfill\qed

\begin{figure}[H]
	\centering
	\includegraphics[width=0.6\textwidth]{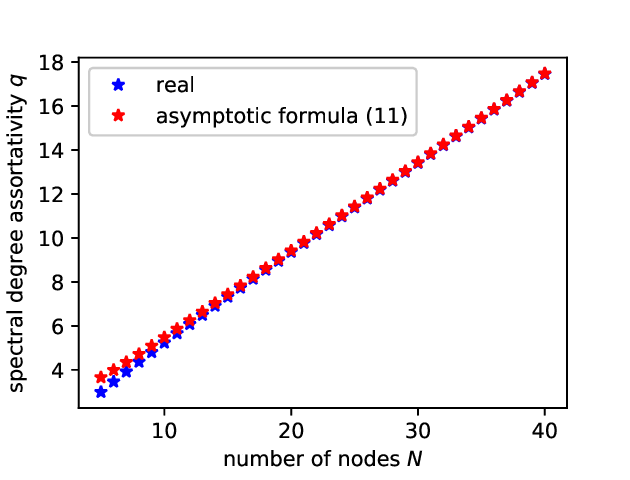}
	\caption{The SDE $q$ for the path graph $P_N$ is computed by simulations (blue) and the approximate relation \eqref{eq_q_path} (red).}
	\label{fig_path_graph_q}
\end{figure}


\subsection{The wheel graph}
\begin{figure}[H]
	\centering
	\includegraphics{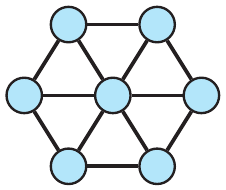}
	\caption{The wheel graph $W_N$. In this example, $N=7$.}
	\label{fig_wheel_graph}
\end{figure}
\begin{theorem}\label{thm_wheel}
	For the \textbf{wheel graph} $W_N$ on $N$ nodes, the SDE $q = 2 + o(1)$ for $N\to\infty$. 
\end{theorem}
\textit{Proof.} The spectral radius of the wheel graph $W_N$ equals (see \cite[p.\ 212]{GS2})
\begin{equation*}
	\lambda_1 = 1 + \sqrt{N}
\end{equation*}
Using that $d_1 = \ldots = d_{N-1} = 3, d_N = N-1$, we find
\begin{equation*}
	(1+\sqrt{N})^q = \frac{N-1}{N} 3^q + \frac{1}{N} (N-1)^q
\end{equation*}
Suppose $q$ converges to some constant $q^*$ for $N\to\infty$, then
\begin{equation*}
	(1+\sqrt{N})^{q^*} = \left( 1 - \frac{1}{N} \right) 3^{q^*} + \frac{1}{N} (N-1)^{q^*}
\end{equation*}
For $N$ large, we may approximate as
\begin{equation*}
	N^{q^*/2} = 3^{q^*} + N^{q^*-1}
\end{equation*}
which can only hold if $q^*/2 = q^*-1$, which results in $q^*=2$. \hfill\qed

\subsection{The "path with double fork" graph}

\begin{figure}[H]
	\centering
	\includegraphics{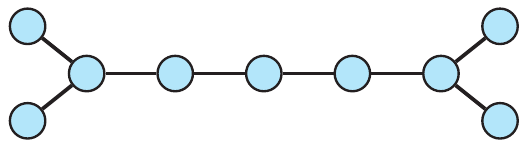}
	\caption{The ``path with double fork graph'' $A_N$ on $N+4$ nodes. In this example, $N=5$.}
	\label{fig_graph_A}
\end{figure}
\begin{theorem}\label{thm_graph_An}
	Define the ``path with double fork" graph $A_N$ on $N+4$ nodes, consisting of a path graph $G$ on $N$ nodes, with two nodes connected to the end of the path and another two nodes attached to the other end of the path. We find that $q \approx 2.36864$, independent of $N$.
\end{theorem}
\textit{Proof.} The spectral radius equals (see \cite[Lemma 3]{cvetkovic1974spectra})
\begin{equation*}
	\lambda_1(A_N) = 2
\end{equation*}
The degree sequence $d_1 = \ldots = d_4 = 1$ and $d_5 = d_6 = 3$ and $d_7 = \ldots = d_{N+4} = 2$. We obtain
\begin{equation*}
	2^q = \frac{4}{N+4} + \frac{N-2}{N+4} 2^{\, q} + \frac{2}{N+4} 3^{\, q}
\end{equation*}
Taking the $2^q$ to the other side yields
\begin{equation*}
	\frac{6}{N+4} 2^{\, q} = \frac{4}{N+4} + \frac{2}{N+4} 3^{\, q}
\end{equation*}
Multiplying by $\frac{N+4}{2}$ finally yields
\begin{equation*}
	3 \cdot 2^{\, q} = 2 + 3^{\, q}
\end{equation*}
This equation has no analytically known solution, but is approximately given by $q \approx 2.36864$. \hfill\qed

The graph $A_N$ is interesting in the sense that for any $N$, both the maximum degree and the spectral radius are fixed. Naturally, the SDE $q$ is then also a constant (independent of $N$).

\subsection{The modified lollipop graph}

\begin{figure}[!ht]
	\centering
	\includegraphics{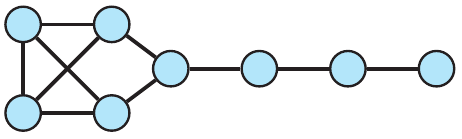}
	\caption{The modified lollipop graph $B_N$ on $N+5$ nodes. Nodes are added to the path on the right-hand side. In this example, $N=3$.}
	\label{fig_graph_B}
\end{figure}
\begin{theorem}
	The modified lollipop graph $B_N$ on $N+5$ nodes is constructed by starting with a complete graph $K_4$ with one link removed. The nodes adjacent to the removed link are connected to a new node~5. To node 5, a path graph of length $N$ is connected, see Figure~\ref{fig_graph_B}. For the graph $B_N$ on $N+5$ nodes, it holds that
	\begin{equation}\label{eq_q_medal}
		q \approx 30.11 \log(N) - 48.46 \qquad \text{for } N \to \infty
	\end{equation}
\end{theorem}
\textit{Proof.} No analytic solution for the spectral radius $\lambda_1$ of $B_N$ seems to exist. However, it appears from simulations that $\lambda_1 \approx 2.902$ for $N \to \infty$ and that convergence to this value is extremely fast. Thus we assume that $\lambda_1$ is a fixed number for all $N$. Substituting the properties of the graph $B_N$ into Eq.~\eqref{eq_def}, we find
\begin{equation*}
	(N+5) \lambda_1^{\, q} = 5 \cdot 3^q + (N-1) 2^q + 1.
\end{equation*}
Assume that $q = a \log(N) + b$, where $a \geq 0$ and $b$ any number. Then we find for large $N$;
\begin{equation*}
	(N+5) \cdot \lambda_1^{a \log(N) + b} = 5 \cdot 3^{a \log(N) + b} + (N-1) 2^{a \log(N) + b} + 1.
\end{equation*}
Using the identity
\begin{equation*}
	c^{a \log(N)} = N^{a \log(c)}
\end{equation*}
for any $a,c,N \geq 0$, we find that
\begin{equation*}
	(N+5) \cdot \lambda_1^b \cdot N^{a \log(\lambda_1)} = 5 \cdot 3^b \cdot N^{a \log(3)} + (N-1) \cdot 2^b \cdot N^{a \log(2)} + 1.
\end{equation*}
In the limit $N\to\infty$, we can approximate $N+5 \approx N$ and $N-1 \approx N$, such that
\begin{equation*}
	\lambda_1^b \cdot N^{1+a \log(\lambda_1)} = 5 \cdot 3^b \cdot N^{a \log(3)} + 2^b \cdot N^{1+a \log(2)} + 1.
\end{equation*}
The second term on the right-hand side is always dominated by the first term on the left-hand side. Also the third term on the right-hand side is dominated by all other terms. Thus, we find
\begin{equation*}
	\lambda_1^b \cdot N^{1+a \log(\lambda_1)} = 5 \cdot 3^b \cdot N^{a \log(3)}
\end{equation*}
Equating the powers of $N$ yields
\begin{equation*}
	1 + a \log(\lambda_1) = a \log(3)
\end{equation*}
which becomes
\begin{equation*}
	a = \frac{1}{\log(3)-\log(\lambda_1)} \approx 30.11
\end{equation*}
Equating the constants gives
\begin{equation*}
	\lambda_1^b = 5 \cdot 3^b
\end{equation*}
which becomes
\begin{equation*}
	b = \frac{\log(5)}{\log\left(\frac{\lambda_1}{3}\right)} \approx -48.46
\end{equation*} 
\hfill\qed

\begin{figure}[H]
	\centering
	\includegraphics[width=0.6\textwidth]{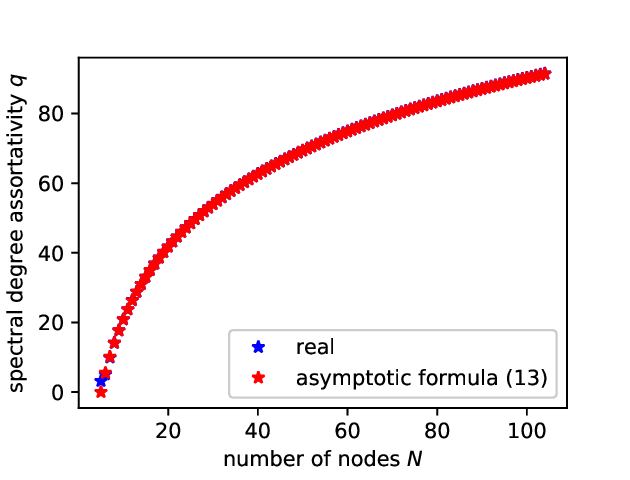}
	\caption{The SDE $q$ for the $B_N$ graph is computed by simulations (blue) and the approximate relation \eqref{eq_q_medal} (red).}
	\label{fig_medal_graph_q}
\end{figure}

Figure~\ref{fig_medal_graph_q} compares the SDE $q$ of the $B_N$ graph obtained using the approximation \eqref{eq_q_medal} (red) and the simulations (blue). The values coincide nearly perfectly, demonstrating the accuracy of the approximation~\eqref{eq_q_medal}. 

The graph $B_N$ is special in the sense that out of all simple, connected graphs with $N=7$ and $N=8$ nodes, the graph $B_N$ has the largest SDE $q$. Provided that $\frac{\lambda_{1}}{d_{\max}}$ and $\frac{d_2}{d_{\max}}$ do not change with the size~$N$, which are approximately satisfied in the $B_N$ graph, we observe that $q$ increases logarithmically in $N$. If the spectral radius $\lambda_{1}$ is close to the maximum degree and $\frac{\lambda_{1}}{d_{\max}}$ decreases with $N$, then $\log\left(  \frac{\lambda_{1}}{d_{\max}}\right)$ is small and $q$ may increase faster than logarithmically in $N$. For example, the path graph was shown to scale $q \sim N$, whereas the $B_N$ graph only scales as $q \sim \log(N)$ for $N \to \infty$.

\section{Simulation results}\label{sec_simulations}
\begin{figure}[H]
	\centering
	\subfloat[\label{fig_11} All graphs $N=7$]{%
		\includegraphics[width=0.48\textwidth]{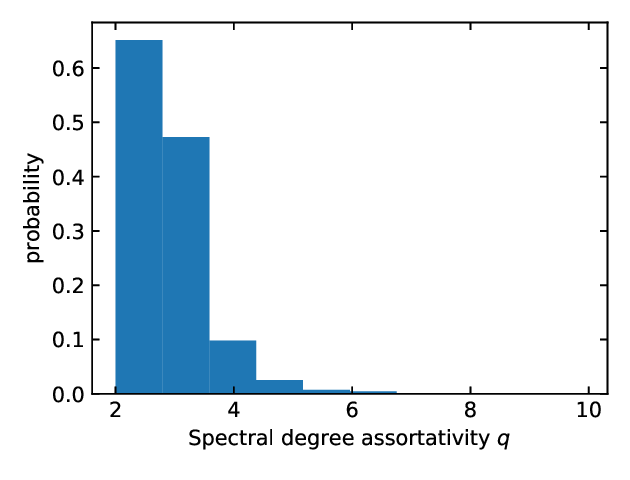}
	}
	\subfloat[\label{fig_12} All graphs $N=8$]{%
		\includegraphics[width=0.48\textwidth]{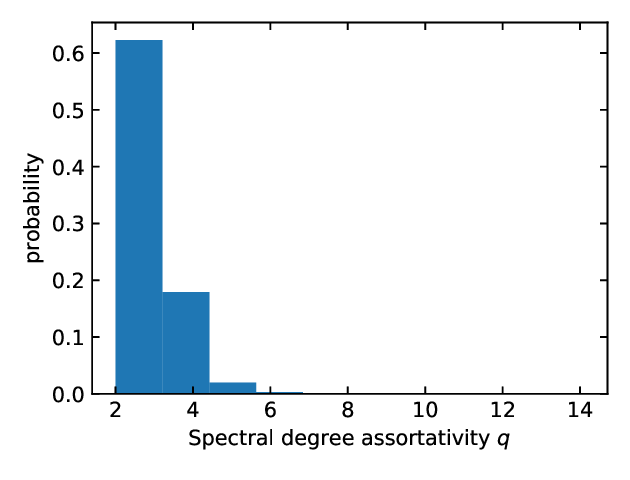}
	} \\
	\subfloat[\label{fig_13} All graphs $N=9$]{%
		\includegraphics[width=0.48\textwidth]{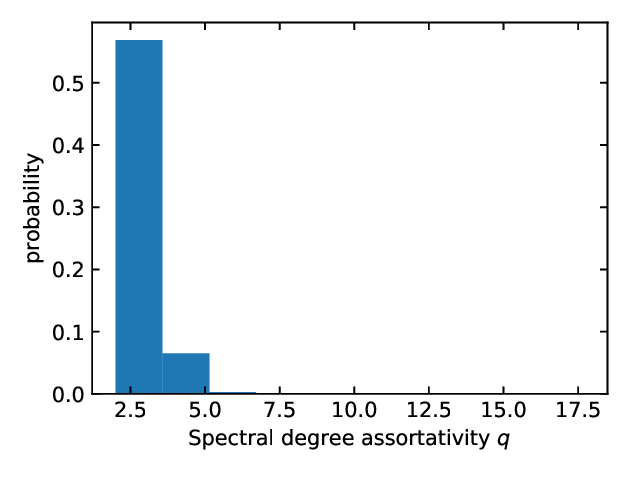}
	}
	\subfloat[\label{fig_14} Topology zoo]{%
		\includegraphics[width=0.48\textwidth]{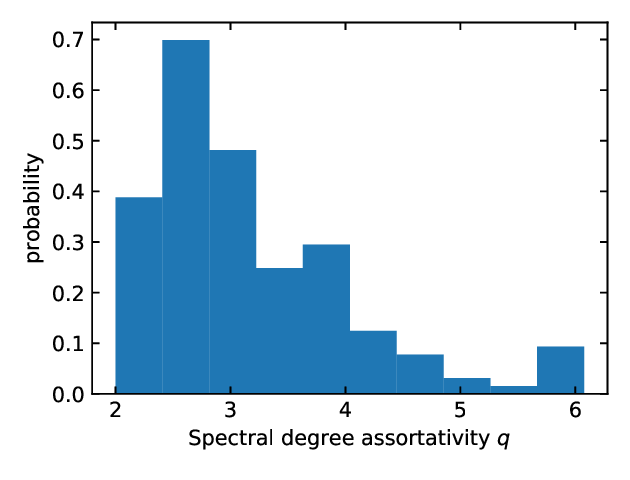}
	} \\
	\subfloat[\label{fig_15} \ER{} graph]{%
		\includegraphics[width=0.48\textwidth]{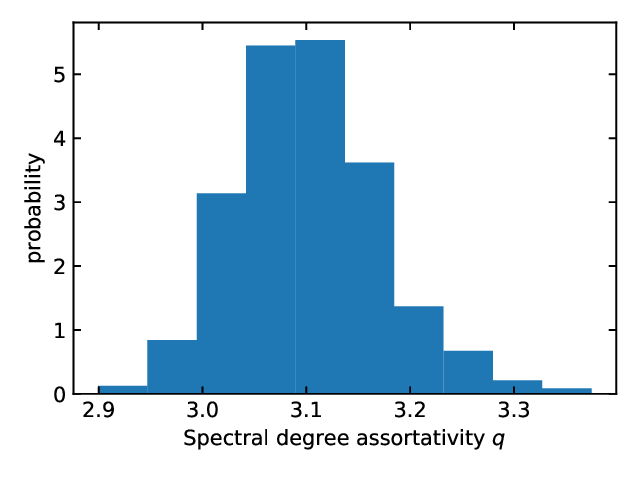}
	}
	\subfloat[\label{fig_16} \BA{} graph]{%
		\includegraphics[width=0.48\textwidth]{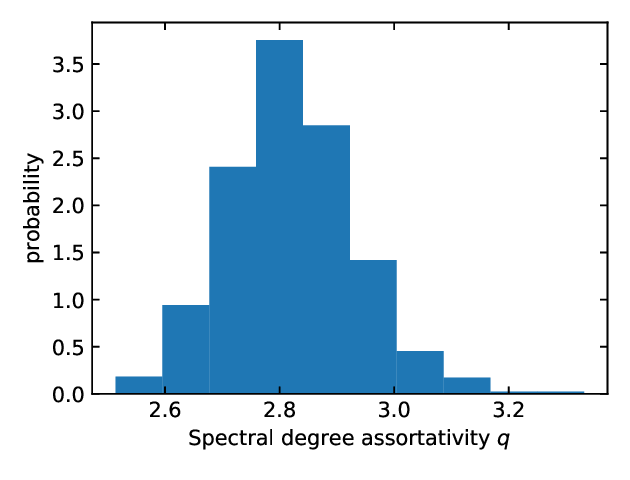}
	}
	\caption{Distribution of the SDE~$q$ for six different graph families. A detailed description of the graph families can be found in Section~\ref{sec_simulations}.}
	\label{fig_hist_q}
\end{figure}

Our main focus in this paper lies on unravelling the interpretation of the SDE $q$ as a graph metric. One approach is to investigate the correlation between the SDE and other graph metrics. We perform an exhaustive search over all non-isomorphic, connected graphs with $N=7, 8$ and $9$ nodes, which are generated using Nauty and Traces \cite{nautytraces}. The number of non-isomorphic, non-regular, connected graphs on $N=7,8,9$ nodes is equal to 849, 11,100 and 261,058, respectively. Additionally, we investigate the real-world communication networks from the Topology Zoo \cite{TopologyZoo}. Out of 261 graphs, we keep the 158 simple, undirected, non-regular, connected graphs. The graphs vary between 4 and 158 nodes with average degree that is at most 4.5. We also generate 1,000 ER graphs on $N=100$ nodes with link probability $p = 0.1 \approx 2 p_c$ and 1,000 \BA{} graphs on $N=100$ nodes with $m=3$. The distribution of the SDE $q$ for all graph families is demonstrated in Figure~\ref{fig_hist_q}. The value of the SDE $q$ varies significantly across the six graph families, but for the most graphs, the value for SDE $q$ is less than four.

\begin{figure}[H]
	\centering
	\subfloat[\label{fig_1} Maximum degree]{%
		\includegraphics[width=0.48\textwidth]{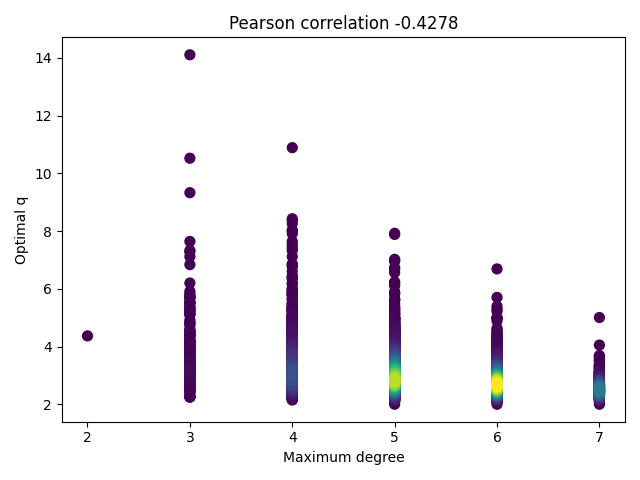}
	}
	\subfloat[\label{fig_2} Degree variance]{%
		\includegraphics[width=0.48\textwidth]{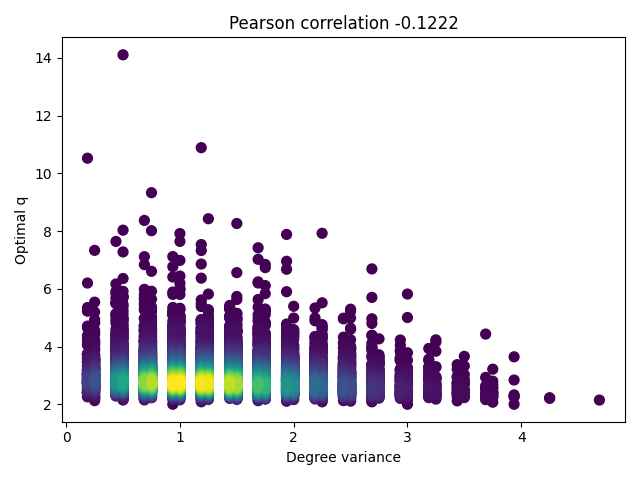}
	} \\
	\subfloat[\label{fig_3} Spectral gap $\lambda_1 -E(D) $]{%
		\includegraphics[width=0.48\textwidth]{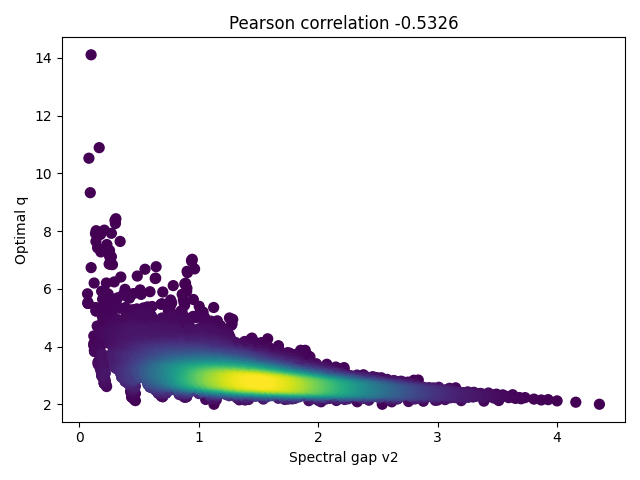}
	}
	\subfloat[\label{fig_4} Algebraic connectivity $\mu_{N-1}$]{%
		\includegraphics[width=0.48\textwidth]{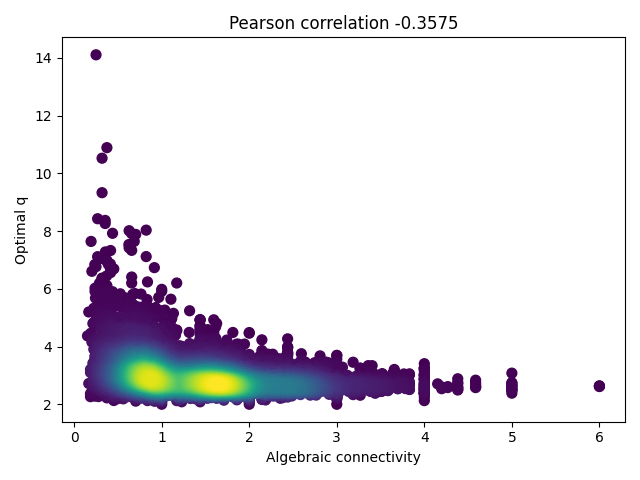}
	} \\
	\subfloat[\label{fig_5} Effective graph resistance $R_G$]{%
		\includegraphics[width=0.48\textwidth]{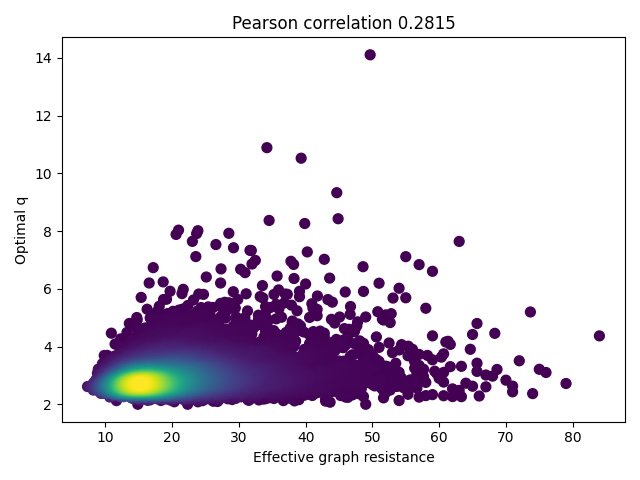}
	}
	\subfloat[\label{fig_6} Degree assortativity]{%
		\includegraphics[width=0.48\textwidth]{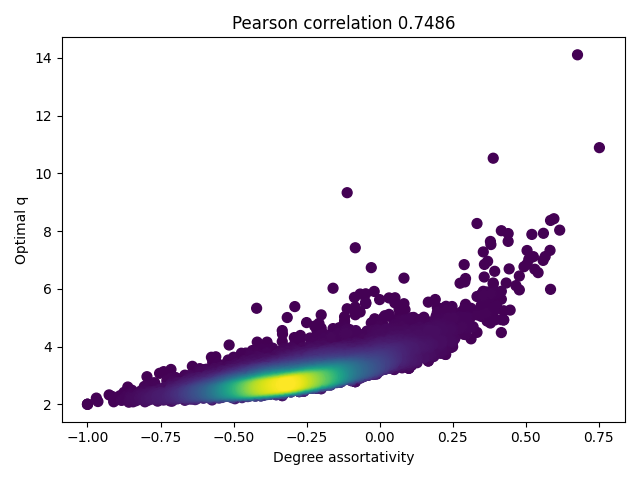}
	}
	\caption{Correlation of SDE $q$ with several graph metrics in all 11,100 connected, non-isomorphic, non-regular graphs with $N=8$ nodes.}
	\label{fig_correlation}
\end{figure}

Table~\ref{tab_correlation} lists the Pearson correlation coefficient for the SDE $q$ with several graph metrics for all graph families. The largest correlation of the SDE $q$ occurs for the spectral gap $\lambda_1 - E[D]$ and the degree assortativity. It was already postulated in Section~\ref{sec_extremal_values} that large values for $q$ arise when $\lambda_1 \approx d_{\max}$, implicitly implying that $q$ will be small if $\lambda_1\ll d_{\max}$. Figure~\ref{fig_correlation} shows the correlation of the SDE~$q$ with a selection of the graph metrics shown in Table~\ref{tab_correlation}. The spectral gap $\lambda_1 - E[D]$ indeed shows high correlation, although not as high as expected.

Most importantly, Table~\ref{tab_correlation} shows that the degree assortativity consistently outperforms the other graph metrics in terms of the Pearson correlation with the SDE. Only for \BA{} graphs, the correlation of the SDE $q$ is stronger with other metrics, which is primarily due to the disassortative nature of \BA{} graphs. All in all, Table~\ref{tab_correlation} confirms the strong correlation between the degree assortativity and the SDE.

\begin{table}[H]
	\centering
	\caption{Correlation between SDE $q$ and several other graph properties, measured in terms of the Pearson correlation coefficient. The correlations are obtained by computing the properties of all non-isomorphic, non-regular connected graphs with $N=7$ (849 graphs), $N=8$ (11,100 graphs) and $N=9$ nodes (261,058 graphs), the Topology Zoo (158 graphs), ER graphs (1,000 random graphs) and BA graphs (1,000 random graphs). Correlations larger than 0.5 (in absolute value) are highlighted in \textbf{bold}.}
	\begin{tabular}{c|c|c|c|c|c|c}
		Metric & $N=7$ & $N=8$ & $N=9$ & Topology Zoo & ER & BA \\
		\hline\hline
		number of links $L$ & -0.100 & -0.144 & -0.179 & 0.311 & -0.140 & 0 \\
		maximum degree $d_{\max}$ & -0.417 & -0.428 & -0.441 & -0.377 & -0.314 & \textbf{-0.773} \\
		minimum degree $d_{\min}$ & -0.294 & -0.323 & -0.364 & 0.234 & -0.197 & -0.007 \\
		degree variance $\textrm{Var}[D]$ & -0.159 & -0.122 & -0.087 & -0.422 & -0.051 & \textbf{-0.855} \\
		spectral radius $\lambda_1$ & -0.068 & -0.100 & -0.127 & -0.196 & -0.058 & \textbf{-0.686} \\
		spectral gap $\lambda_1-\lambda_2$ & -0.135 & -0.131 & -0.133 & -0.447 & -0.021 & -0.423 \\
		spectral gap $\lambda_1 - E[D]$ & \textbf{-0.535} & \textbf{-0.533} & \textbf{-0.528} & -0.394 & -0.331 & \textbf{-0.765} \\
		algebraic connectivity $\mu_{N-1}$ & -0.318 & -0.357 & -0.397 & -0.321 & -0.243 & -0.081 \\
		effective graph resistance $R_G$ & 0.221 & 0.282 & 0.324 & 0.309 & 0.228 & -0.405 \\
		average shortest path length & 0.280 & 0.327 & 0.356 & \textbf{0.606} & 0.231 & \textbf{0.892} \\
		diameter & 0.422 & 0.437 & 0.444 & \textbf{0.635} & 0.142 & 0.240 \\
		clustering coefficient & 0.167 & 0.131 & 0.093 & 0.145 & 0.053 & -0.123 \\
		radius (min eccentricity) & 0.338 & 0.303 & 0.272 & \textbf{0.625} & 0 & 0.089 \\
		size of largest clique & 0.158 & 0.128 & 0.105 & 0.131 & 0.037 & -0.078 \\
		node connectivity & -0.268 & -0.313 & -0.360 & 0.207 & -0.197 & -0.007 \\
		link connectivity & -0.300 & -0.330 & -0.368 & 0.234 & -0.197 & -0.007 \\
		degree assortativity & \textbf{0.765} & \textbf{0.749} & \textbf{0.747} & \textbf{0.788} & \textbf{0.856} & \textbf{0.712} \\
		number of bridges & 0.206 & 0.286 & 0.345 & -0.281 & 0.112 & 0.001 \\
		local efficiency & -0.104 & -0.166 & -0.223 & 0.063 & -0.065 & \textbf{-0.625} \\
		global efficiency & -0.183 & -0.231 & -0.263 & -0.459 & -0.208 & \textbf{-0.906} \\
		number of leaf nodes & 0.128 & 0.226 & 0.307 & -0.341 & 0.112 & 0.001 \\
		graph energy & 0.071 & -0.044 & -0.127 & 0.479 & -0.125 & \textbf{0.662} \\
		Estrada index & -0.074 & -0.086 & -0.104 & 0.001 & -0.072 & \textbf{-0.656} \\
		number of spanning trees & -0.122 & 0.116 & -0.121 & 0.120 & -0.055 & \textbf{0.535} \\
		max. Laplacian eigenvalue & -0.398 & -0.417 & -0.436 & -0.356 & -0.294 & \textbf{-0.779} \\
		\hline\hline
	\end{tabular}
	\label{tab_correlation}
\end{table}

\section{Conclusion}
\label{sec_conclusion}
In this paper, we introduced the spectral degree exponent (SDE) as a novel graph metric. We guaranteed the existence and uniqueness of the metric. Moreover, the extremal values for both metrics are attained in very similar graphs. We further proceeded by deriving bounds on and iterative formulas for the SDE. We also (asymptotically) computed the SDE for several graph families. Finally, we showed that amongst many known metrics, the degree assortativity shows the largest correlation with the degree assortativity in a plethora of graphs, ranging from small connected graphs, real-world communication networks and \ER{} and \BA{} graphs.

We see several directions for future research. Since the spectral degree exponent is defined implicitly, tight upper and lower bounds are welcome. Several other properties of the SDE have not yet been established, e.g.\ whether the SDE is submodular or if the SDE can be used successfully in applications to enhance or reduce assortativity. We also welcome exact or asymptotic solutions of the SDE on other graph families.

\section*{Acknowledgments}

This research has been funded by the European Research Council (ERC) under the European Union's Horizon 2020 research and innovation programme (grant agreement No 101019718).


\bibliographystyle{comnet}

\begin{thebibliography}{00}
	
	\bibitem{BrouwerHaemers}
	Brouwer, A.~E. {\&} Haemers, W.~H. (2011) {\em {Spectra of graphs}}.
	Springer Nature, New York, United States, first edition.
	
	\bibitem{cvetkovic1974spectra}
	Cvetkovi{\'c}, D., Gutman, I. {\&} Trinajsti{\'c}, N. (1974)  Conjugated
	molecules having integral graph spectra. {\em Chemical Physics Letters},
	\textbf{29}(1), 65--68.
	
	\bibitem{HofmeisterThesis}
	Hofmeister, M. (1986) {\em {\"U}ber {Z}usammenh{\"a}nge zwischen
		{V}alenzsequenz, {E}ckenzahl und {G}r{\"o}{\ss}ten {E}igenwerten von
		{G}raphen und {D}igraphen}.
	Phd., Universit{\"a}t zu K{\"o}ln.
	
	\bibitem{hofmeister1988spectralradius}
	Hofmeister, M. (1988)  {Spectral Radius and Degree Sequence}. {\em
		Mathematische Nachrichten}, \textbf{139}(1), 37--44.
	
	\bibitem{TopologyZoo}
	Knight, S., Nguyen, H.~X., Falkner, N., Bowden, R. {\&} Roughan, M. (2011)
	{The Internet Topology Zoo}. {\em Selected Areas in Communications, IEEE
		Journal on}, \textbf{29}(9), 1765--1775.
	
	\bibitem{nautytraces}
	McKay, B.~D. {\&} Piperno, A. (2014)  {Practical graph isomorphism, II}. {\em
		Journal of Symbolic Computation}, \textbf{60}, 94--112.
	
	\bibitem{newman2002assortativity}
	Newman, M. E.~J. (2002)  {Assortative Mixing in Networks}. {\em Phys. Rev.
		Lett.}, \textbf{89}, 208701.
	
	\bibitem{BookNewman}
	Newman, M. E.~J. (2010) {\em {Networks: An Introduction}}.
	Oxford University Press, Oxford, United Kingdom, first edition.
	
	\bibitem{noldus2015assortativity}
	Noldus, R. {\&} Van~Mieghem, P. (2015)  Assortativity in complex networks. {\em
		Journal of Complex Networks}, \textbf{3}(4), 507--542.
	
	\bibitem{GS2}
	Van~Mieghem, P. (2023) {\em {Graph Spectra for Complex Networks}}.
	Cambridge University Press, Cambridge, United Kingdom, second edition.
	
	\bibitem{vanmieghem2010degree}
	Van~Mieghem, P., Wang, H., Ge, X., Tang, S. {\&} Kuipers, F.~A. (2010)
	Influence of assortativity and degree-preserving rewiring on the spectra of
	networks. {\em The European Physics Journal B}, \textbf{76}, 643--652.
	
\end{thebibliography}

%


\appendix

\section{Proof of Theorem~\ref{thm_q=2}}\label{app_proof_q2}
We follow here the proof of Hofmeister \cite{hofmeister1988spectralradius}, but extend it for weighted graphs.

(a) If $G$ is biregular, its spectral radius $\lambda_1 = \sqrt{r_1 r_2}$, where $r_1$ and $r_2$ are the two weighted degrees \cite{hofmeister1988spectralradius}. If the size of the groups are $m$ and $n$ respectively, we recover the following relation; $r_2 = r_1 \frac{m}{n}$ due to the biregular constraint. Using $m+n=N$, we deduce that
\begin{align*}
	\frac{n}{N} = \frac{r_1}{r_1 + r_2}, \qquad \frac{m}{N} = \frac{r_2}{r_1 + r_2}.
\end{align*}
Since $d_1 = \ldots d_m = r_1$ and $d_{m+1} = \ldots = d_N = r_2$ for a biregular graph, we find from \eqref{def_q} for $q=2$ that
\begin{align*}
	\lambda_1 &= \left( \frac{1}{N} \sum_{i=1}^N d_i^2 \right)^{1/2} = \left( \frac{m}{N} r_1^2 + \frac{n}{N} r_2^2 \right)^{1/2} = \left( \frac{r_1 ^2 r_2}{r_1 + r_2} + \frac{r_1 r_2^2}{r_1 + r_2} \right)^{1/2} = \sqrt{r_1 r_2}
\end{align*}
which proves that the SDE $q=2$ for biregular graphs.

(b) Suppose that $G$ is connected and $q=2$. We first prove that the graph $G$ contains exactly two weighted degrees.

Consider the following optimization problem.
\begin{equation}\label{eq_optimisation_problem}
	\begin{split}
		\max \quad &\| A \x \|^2 \\
		\text{s.t.} \quad &\|\x\|^2 = 1
	\end{split}
\end{equation}
The optimum is attained when $\x$ equals the principal eigenvector of the adjacency matrix~$A$. Choosing $\x_0 = \frac{1}{\sqrt{N}} \u$, we find
\begin{equation*}
	\| A \x_0 \|^2 = \frac{1}{N} (A \u)^T (A \u) = \frac{1}{N} \sum_{i=1}^N d_i^2 = \lambda_1^2
\end{equation*}
where the latter equality follows due to $q=2$. Thus, $\x_0$ is indeed the principal eigenvector of $A$.

We can reformulate the optimization problem \eqref{eq_optimisation_problem} as a Lagrange optimization problem. Define $\Phi(\x) = \| A \x\|^2$ and $\Psi(\x)=\|\x\|^2$. Then there must exist a Lagrange multiplier $\xi$ such that the Lagrangian
\begin{equation*}
	\mathcal{L}(\x, \xi) = \Phi(\x) - \xi \Psi(\x)
\end{equation*}
is minimal at $\x_0$. In other words,
\begin{equation*}
	\pd{\mathcal{L}}{x_k}(\x_0) = 0.
\end{equation*}
Using
\begin{equation*}
	\pd{\Phi(\x)}{x_k} = 2 \sum_{i=1}^N a_{ik} \sum_{j=1}^N a_{ij} x_j, \qquad \pd{\Psi(\x)}{x_k} = 2 x_k, \qquad \x_0 = \frac{1}{\sqrt{N}} \u,
\end{equation*}
we find
\begin{equation*}
	\xi = \sum_{i=1}^N a_{ik} d_i
\end{equation*}
In terms of unweighted graphs, the Lagrange multiplier $\xi$ is equal to the number of two-hop walks starting at node $k$.

Since $G$ is not regular, there are at least two different weighted degrees: the maximal weighted degree $d_{\max}$ and minimum weighted degree $d_{\min}$. Suppose there are more than two weighted degrees, so there exists a node $u$ with weighted degree $d_{\min} < d_u < d_{\max}$. Since $G$ is connected, node $u$ must be connected to a node of weighted degree $d_{\min}$ or $d_{\max}$.

Suppose node $u$ is connected to a node $v$ of weighted degree $d_{\max}$. Then
\begin{align*}
	\xi &= \sum_{i=1}^N a_{iv} d_i \\
	&= a_{uv} d_u + \sum_{i=1, i \neq u}^N a_{iv} d_i \\
	&\geq a_{uv} d_u + \sum_{i=1, i \neq u}^N a_{iv} d_{\min} \\
	&= a_{uv} d_u - a_{uv} d_{\min} + d_{\min} \sum_{i=1}^N a_{iv} \\
	&= a_{uv} (d_u - d_{\min}) + d_{\min} d_{\max} \\
	&> d_{\min} d_{\max}
\end{align*}
Conversely, suppose node $u$ is connected to a node $v$ of weighted degree $d_{\min}$. Then
\begin{align*}
	\xi &= \sum_{i=1}^N a_{iv} d_i \\
	&\leq \sum_{i=1}^N a_{iv} d_{\max} \\
	&= d_{\min} d_{\max}
\end{align*}
Using the Lagrange multiplier $\xi$, we conclude that the graph $G$ cannot contain more two different weighted degrees. Hence, the graph $G$ contains two weighted degrees $r_1 < r_2$, as exemplified by the Lagrange multiplier;
\begin{align*}
	\xi &= \sum_{i=1}^N a_{iu} d_i \geq r_2 \sum_{i=1}^N a_{iu} = r_1 r_2
\end{align*}
and
\begin{align*}
	\xi &= \sum_{i=1}^N a_{iv} d_i \leq r_1 \sum_{i=1}^N a_{iv} = r_1 r_2
\end{align*}
It remains to show that nodes with weighted degree $r_1$ and $r_2$ are not connected to each other. Suppose node $u$ and $v$ are adjacent and both have weighted degree $r_1$. Then 
\begin{equation*}
	\xi = \sum_{i=1}^N a_{ui} d_i < \sum_{i=1}^N a_{ui} r_2 = r_1 r_2
\end{equation*}
contradicting the earlier obtained result that $\xi = r_1 r_2$. In a similar fashion, one can show that nodes with degrees $r_2$ cannot be connected to each other. \hfill $\square$

\section{Proof of Theorem~\ref{thm_q=infty}}\label{app_proof_qinfty}
Suppose $q=\infty$. As shown in Eq.\ \eqref{eq_relation_dmax}, if $q\to\infty$, then $\lambda_1 = d_{\max}$. All graphs that satisfy $\lambda_1 = d_{\max}$ must have a component which is a clique of size $d_{\max}+1$. From all connected graphs, only the weighted regular graphs satisfy this property, but the regular graphs have been explicitly excluded in the theorem statement. Thus we conclude $G$ is disconnected.

Now suppose the reverse holds, i.e.\ $G$ contains a clique of weight $d_{\max}+1$. Again, $G$ cannot be a complete graph, so $G$ must be disconnected. Since $G$ consists of multiple components, the spectral radius $\lambda_1$ of $G$ is the maximum spectral radius of its components. Since one component is a clique of size $d_{\max}+1$ and the maximum weighted degree is $d_{\max}$, this component inevitably has the largest spectral radius. Thus the spectral radius of $G$ equals $d_{\max}$. Substituting $\lambda_1 = d_{\max}$ into \eqref{eq_relation_dmax} yields
\begin{equation*}
	d_{\max} = \left( \frac{c}{N} d_{\max}^{\, q} + \frac{1}{N} \sum_{i=c+1}^N d_i^{\, q} \right)^{1/q}
\end{equation*}
Rearranging yields
\begin{equation*}
	d_{\max}^{\, q} = \frac{c}{N} d_{\max}^{\, q} + \frac{1}{N} \sum_{i=c+1}^N d_i^{\, q}
\end{equation*}
and can be further simplified to
\begin{equation*}
	N - c = \sum_{i=c+1}^N \left( \frac{d_i}{d_{\max}} \right)^{\, q}
\end{equation*}
Since $d_i < d_{\max}$ any finite value of $q$ yields $\left( \frac{d_i}{d_{\max}} \right)^{\, q} < 1$. Since the sum is composed of $N-c$ terms, no finite $q$ can yield equality. Thus we conclude that $q=\infty$. \hfill\qed

\end{document}